\documentclass[12pt]{article}

\usepackage{amsmath,amssymb,amsfonts,amsthm,epic}

\usepackage{graphicx}

\newcommand\torder{\vartriangleright}

\theoremstyle{plain}
\newtheorem{thm}[subsubsection]{Theorem}

\newtheorem{ex}[subsubsection]{Example}

\theoremstyle{definition}
\newtheorem{rem}[subsubsection]{Remark}
\newtheorem{defn}[subsubsection]{Definition}
\newtheorem{nota}[subsubsection]{Notation}

\newcommand\brem{\begin{rem}\begin{sffamily}\begin{upshape}}
\newcommand\erem{\end{upshape}\end{sffamily}\end{rem}}
\newcommand\bdefn{\begin{defn}\begin{rm}}
\newcommand\edefn{\end{rm}\hfill$\Box$\end{defn}}
\newcommand\bnot{\begin{nota}\begin{rm}}
\newcommand\enot{\end{rm}\hfill$\Box$\end{nota}}
\newcommand\bex{\begin{ex}\begin{rm}}
\newcommand\eex{\end{rm}\hfill$\Box$\end{ex}}

\newenvironment{Proof}{%
\par\noindent{\scshape Proof:}\begin{rm}}{\hfill$\Box$\end{rm}\newline}
\numberwithin{equation}{subsection}
\title{Irrational HK multiplicities are possible for trinomial hyper surfaces}
\date{}
\author {Shyamashree Upadhyay\\Department of Mathematics\\ Indian Institute of Technology, Guwahati\\Assam-781039, INDIA\\email: shyamashree@iitg.ernet.in}
\begin{document}
\maketitle
\begin{abstract}
A `trinomial hyper surface' is defined in \S \ref{s.Introduction} below. In this article, I provide a supportive reasoning towards the fact that there can be examples of trinomial hyper surfaces (at least over fields of characteristic $2$) for which the corresponding Hilbert-Kunz multiplicity can become irrational. 
\end{abstract}
\tableofcontents
\section{Introduction}\label{s.Introduction}
Let $S=k[x_1,\ldots,x_m]$ where $k$ is a field of arbitrary prime characteristic $p>0$ and $J=(f)$ be an ideal in $S$ where $f$ is a polynomial containing $3$ non-constant terms in it and no constant term. The affine variety defined by the ideal $J$ is called a \textit{trinomial hyper surface}. For a detailed definition, see \S 2.0.1 of \cite{Su1}. 

In the present work, I explain how there can exist examples of trinomial hyper surfaces (at least over fields of characteristic $2$) for which the corresponding Hilbert-Kunz multiplicity can become irrational. Before this, I have done a similar work regarding computation of the Hilbert-Kunz function for `disjoint-term' trinomial hyper surfaces (see \cite{Su1} and \cite{Su2}). In \cite{Su1} and \cite{Su2}, my strategy of work was to first apply the process `mutation' as mentioned in \cite{Su} to reduce the problem to a problem of solving certain systems of linear equations. Then by applying some basic linear algebra techniques on the resulting problem (of solving certain systems of linear equations), we were reduced to a problem of computing the ranks of a huge collection of matrices (which are of size of the order of $p^n$). This rank computation problem was solved in \cite{Su2}. But for the case of `disjoint-term' trinomial hyper surfaces, it turns out that the corresponding Hilbert-Kunz multiplicity will be rational, see \S 5 of \cite{Su2} for a reasoning. 

My strategy here is to apply a similar (as in \cite{Su1}) technique on any trinomial hyper surface in general (which is not necessarily having `disjoint terms' in it). Then as in the case considered in \cite{Su1} and \cite{Su2}, we will encounter a similar-looking rank computation problem here. But here, the nature of the matrices (which are involved in the rank computation problem) will tell us why there can exist examples of trinomial hyper surfaces (at least over fields of characteristic $2$) for which the corresponding Hilbert-Kunz multiplicity can become irrational.

\noindent\textbf{An apology:} The write-up of this article is not very detailed or precise. The reason behind this is that I wrote it in a hurry. I deeply regret for this.
\section{The main theorem}\label{s.the-main-theorem}
In \cite{Su1}, the problem of computing the Hilbert-Kunz function for disjoint term trinomial hyper surfaces was first translated to a theorem (theorem 3.2.3 of \cite{Su1}). After this translation, the problem was reduced to linear algebra and combinatorics. In this section, I state the main theorem for trinomial hyper surfaces (theorem \ref{t.mutation-stop} below) which is the analog of theorem 3.2.3 of \cite{Su1}. But before stating the main theorem, we need to introduce some terminology, notation e.t.c..
\subsection{The term order $\torder$}\label{ss.torder}
Let us put an order (denote it by $\torder$) on the set of all monomials in the variables $x_1,\ldots,x_m$ as follows:---
\begin{itemize}
\item Set $x_1\torder\cdots\torder x_m$.
\item On the set of all monomials in the variables $x_1,\ldots,x_m$, $\torder$ is the degree lexicographic order with respect to the order $\torder$ defined on the variables $x_1,\ldots,x_m$.
\end{itemize} 
Since the polynomial $f$ has $3$ non-constant terms in it, let us denote the most initial (with respect to $\torder$) term of $f$ as $[3]$, the next most initial term of $f$ as $[2]$ and the least initial term as $[1]$. Hence we have
$$J=(f)=([3]+[2]+[1])$$
We can arrange the variables $x_1,\ldots,x_m$ in such a way that the ordered set $x_1\torder\cdots\torder x_m$ is the same as the ordered set $E_q\torder\cdots\torder E_1\torder N_r\torder\cdots\torder N_1\torder Z_s\torder\cdots\torder Z_1\torder P_1\torder\cdots\torder P_t$ where $N_r,\ldots,N_1,Z_s,\ldots,Z_1,P_1,\ldots,P_t$ are the negative, zero and the positive difference variables respectively with respect to the terms $[2]$ and $[1]$ of $f$ as mentioned in \S 4.2.1 of \cite{Su}. And the variables $E_q,\ldots,E_1$ are the \textit{extra} variables which appear in the term $[3]$ only and not in the other terms $[2]$ and $[1]$. Of course, there can be cases where these extra variables do not exist. But we can assume without loss of generality that they exist. The variables $E_q,\ldots,E_1$ are having the property that the $deg_{f}(E_q)\geq\cdots\geq deg_{f}(E_1)$ where for any $i\in\{1,\ldots,q\}$, $deg_{f}(E_i)$ denotes the degree of the variable $E_i$ in the polynomial $f$. 
\brem\label{r.scalarsintermsoff}
Note here that the terms $[3],[2],[1]$ of $f$ are assumed to be containing scalar coefficients.
\erem 
Recall the set $\mathfrak{M}:=\{a_t|t\in\{1,\ldots,p^{mn}\}\}$ from \S 2.1 of \cite{Su}. Let $A$ be an arbitrary element of the set $\mathfrak{M}$. Recall the ideal $A_c$ corresponding to the element $A$ from \S 3.2 of \cite{Su1}. The definitions of the set $\mathfrak{M}$ and the ideal $A_c$ corresponding to any element $A$ of it are similar for the present case.  

Given any term $\tau$ of the polynomial $f$, define $[-\tau]:=\frac{1}{\tau}$. Let
\begin{center}
$\mathcal{B}(A,f):=$ the set of all non-convergent mutants in $A$ and $f$ (see \S 3.2 of \cite{Su} for the definition of a mutant)\\
and $\mathcal{A}(A,f):=\{B[-\tau]|B\in\mathcal{B}(A,f)\ and\ \tau\ divides\ B\}$. 
\end{center}
As in \cite{Su1}, the problem of computing the Hilbert-Kunz function reduces to determining whether $A\in A_c+J$ or not for each monomial $A\in\mathfrak{M}$. Theorem \ref{t.mutation-stop} below provides the answer to this:
\begin{thm}\label{t.mutation-stop}
Let $f=[3]+[2]+[1]$. $A\in A_c+J$ iff one of the following mutually exclusive conditions hold:---\\
(i) The term $[1]$ of $f$ divides the monomial $A$.\\
\noindent (ii) The term $[1]$ of $f$ doesnot divide the monomial $A$, but the term $[2]$ divides $A$, there exist extra variables $E_1,\ldots,E_q$ in the term $[3]$ and there exists a positive integer $M$ for which $A\dfrac{{[-2]}^{M-1}}{{[1]}^{M-1}}[-2]$ contains no negative powers and $A\dfrac{{[-2]}^M}{{[1]}^{M}}$ is convergent. \\
\noindent (iii) Neither condition (i) nor condition (ii) holds, but at least one of the terms $[3]$ and $[2]$ divides $A$ and there exists scalars $c_D$ (corresponding to each $D\in\mathcal{A}(A,f)$) such that the product $f.(\Sigma_{D\in\mathcal{A}(A,f)}c_{D}D)$ equals $bA+finitely\ many\ convergent\ terms$ for some non-zero scalar $b$.
\end{thm}
\begin{Proof}
Similar to the proof of theorem 3.2.3 of \cite{Su1}. 
\end{Proof}
\section{Reduction to linear algebra}\label{s.red-lin-alg}
We will now provide an equivalent formulation of condition (iii) of theorem \ref{t.mutation-stop} above in terms of linear algebra. This equivalent formulation follows from remark \ref{r.condition(iii)-to-sysoflineareqn} below.
\brem\label{r.condition(iii)-to-sysoflineareqn}
Let $\mathcal{E}(A,f):=\{D[\tau]|D\in\mathcal{A}(A,f)\ and\ \tau\in\{[1],[2],[3]\}\}$ and $\mathcal{L}(A,f):=$the set of all elements in $\mathcal{E}(A,f)$ which when expressed in lowest terms do not contain all $[3]$ s in the denominator. Note that the product $f.(\Sigma_{D\in\mathcal{A}(A,f)}c_{D}D)$ in condition (iii) of theorem \ref{t.mutation-stop} above equals a linear combination of elements of $\mathcal{E}(A,f)$, say $\Sigma_{B\in\mathcal{E}(A,f)}e_{B}B$. This sum can be broken into $2$ parts as follows: $\Sigma_{B\in\mathcal{E}(A,f)}e_{B}B=\Sigma_{B\notin\mathcal{L}(A,f)}e_{B}B+\Sigma_{B\in\mathcal{L}(A,f)}e_{B}B$. Since any mutant in $A$ and $f$ (expressed in lowest terms) that contains only $[-2]$'s in the numerator and only $[3]$ s in the denominator is convergent, it follows that the portion $\Sigma_{B\notin\mathcal{L}(A,f)}e_{B}B$ contains all convergent terms and there are finitely many such terms. So if we equate the coefficients of like terms of the product $f.(\Sigma_{D\in\mathcal{A}(A,f)}c_{D}D)$ and the sum $\Sigma_{B\in\mathcal{L}(A,f)}e_{B}B$, we get a system $\mathfrak{A}_{A,f}\mathfrak{X}=\mathfrak{B}$ of linear equations where $\mathfrak{X}$ is a column vector in the unknowns $c_D$, $\mathfrak{B}$ is a column vector in the scalars $e_B$ where $B\in\mathcal{L}(A,f)$ and $\mathfrak{A}_{A,f}$ is a matrix with entries from the set $\{0,1\}$. 
\erem
So an equivalent formulation of condition (iii) of theorem \ref{t.mutation-stop} above will be condition  $(iii)'$ as stated  below: 
\begin{quote}
$(iii)'$ Neither condition (i) nor condition (ii) of theorem \ref{t.mutation-stop} holds, but at least one of the terms $[3]$ and $[2]$ divides $A$ and the system $\mathfrak{A}_{A,f}\mathfrak{X}=\mathfrak{B}$ of linear equations is solvable for the vector $\mathfrak{B}$ which is having the property that:
\begin{center}
$e_B\neq 0$ for $B=A$ and $e_B=0$ for all $B$ non-convergent.
\end{center}
\end{quote}
Combining theorem \ref{t.mutation-stop} and condition $(iii)'$ above, we get the following theorem:
\begin{thm}\label{t.mutation-stop-modified}
$A\in A_c+J$ iff 
\begin{itemize}
\item either condition (i) or condition (ii) of theorem \ref{t.mutation-stop} holds or
\item neither condition (i) nor (ii) of theorem \ref{t.mutation-stop} holds but at least one of the terms $[2]$ and $[3]$ divides $A$ and the system $\mathfrak{A}_{A,f}\mathfrak{X}=\mathfrak{B}$ of linear equations is solvable for the vector $\mathfrak{B}$ (this vector is introduced in remark \ref{r.condition(iii)-to-sysoflineareqn} above) which is having the property that:
\end{itemize}
\begin{center}
$e_B\neq 0$ for $B=A$ and $e_B=0$ for all $B$ non-convergent.
\end{center}
\end{thm}
\subsection{Further reduction to another linear system}\label{ss.further-red-toanother}
Let $A$ be a monomial such that neither condition (i) nor condition (ii) of theorem \ref{t.mutation-stop} hold. Consider the system $\mathfrak{A}_{A,f}\mathfrak{X}=\mathfrak{B}$ of linear equations where the indexing set for the column vector $\mathfrak{X}$ is $\mathcal{A}(A,f)$ and that for the column vector $\mathfrak{B}$ is $\mathcal{L}(A,f)$. In the present case of trinomial hyper surfaces, these indexing sets $\mathcal{A}(A,f)$ and $\mathcal{L}(A,f)$ happen to be infinite, unlike in the case of disjoint term trinomial hyper surfaces which was discussed in \cite{Su1}. However, I am not providing a detailed description of these indexing sets here.

For any fixed $B\in \mathcal{L}(A,f)$, look at the non-zero entries in the matrix $\mathfrak{A}_{A,f}$ in the row corresponding to it. These non-zero entries appear in at most $3$ columns of the matrix $\mathfrak{A}_{A,f}$. These columns are arranged in the increasing order induced by the order $\preceq$ on $\mathcal{A}(A,f)$ [see \S 3 of \cite{Su1} for the definition of the order $\preceq$].  For the row corresponding to the fixed $B\in \mathcal{L}(A,f)$, look at all those column(s) in the matrix $\mathfrak{A}_{A,f}$ which have non-zero entries. And out of all these column(s) choose the `smallest´ with respect to the order $\preceq$ on $\mathcal{A}(A,f)$. Denote this column as $C_{1,B}$.

There may exist two distinct elements $B, \acute{B}\in\mathcal{L}(A,f)$ for which the columns $C_{1,B}$ and $C_{1,\acute{B}}$ are the same. 
\begin{center}
Let $\mathcal{R}(A,f)$ denote the set of all elements $B\in\mathcal{L}(A,f)$ for which there exists an element $\acute{B}\in\mathcal{L}(A,f)$ distinct from $B$ such that $C_{1,B}$ matches with $C_{1,\acute{B}}$ and such that $\acute{B}\precsim B$ with respect to the order $\precsim$ on $\mathcal{L}(A,f)$. [see \S 3.3 of \cite{Su1} for a definition of the order $\precsim$.]\\
Let $\mathcal{S}(A,f):=\mathcal{L}(A,f)-\mathcal{R}(A,f)$. \\
The sets $\mathcal{S}(A,f)$ and $\mathcal{R}(A,f)$ are infinite.
\end{center}
Now by applying the same techniques as in \S 3.3 of \cite{Su1}, we transform the system $\mathfrak{A}_{A,f}\mathfrak{X}=\mathfrak{B}$ to another equivalent system $\mathfrak{A}_{A,f}'\mathfrak{X}=\mathfrak{B}'$ of linear equations. And then by applying the same techniques as in \S 4.1 of \cite{Su1}, we reduce the system $\mathfrak{A}_{A,f}'\mathfrak{X}=\mathfrak{B}'$ to another equivalent system $\mathfrak{B}_{A,f}\mathfrak{Y}=\mathfrak{e}$ of linear equations whose description is given in subsection \ref{ss.description-of-anewsys} below.
\section{Reduction to Combinatorics}\label{s.red-to-comb}
The system $\mathfrak{B}_{A,f}\mathfrak{Y}=\mathfrak{e}$ of linear equations obtained at the end of subsection \ref{ss.further-red-toanother} has nice combinatorial properties which will give us a reasoning behind the suspicion that irrationality of HK multiplicity can happen in this case. We therefore need to have an understanding of this system of linear equations in details. This description is given in subsection \ref{ss.description-of-anewsys} below.  
\subsection{Description of the new system}\label{ss.description-of-anewsys}
The new system $\mathfrak{B}_{A,f}\mathfrak{Y}=\mathfrak{e}$ of linear equations can be described as follows:\\
\begin{itemize}
\item The number of rows in the matrix $\mathfrak{B}_{A,f}$ is $|\mathcal{R}(A,f)|$. The rows of the matrix $\mathfrak{B}_{A,f}$ are indexed by elements of the set $\mathcal{R}(A,f)$ which are arranged in the decreasing order induced by the order $\precsim$ on $\mathcal{L}(A,f)$.
\item The number of columns in it is $|P_{A,f}|$ where $P_{A,f}:=$ the set of all convergent elements in $\mathcal{L}(A,f)$. The columns of the matrix $\mathfrak{B}_{A,f}$ are indexed by elements of the set $P_{A,f}$.
\item The column vector $\mathfrak{e}$ is given by $[\ldots,0,0,e_A]^t$ where $e_A\neq 0$. 
\item The entries of the column vector $\mathfrak{Y}$ are the elements of the vector $\mathfrak{B}$ which correspond to the elements of the set $P_{A,f}$ and they are arranged in the increasing order induced by the order $\precsim$ on $\mathcal{L}(A,f)$. 
\item The entries of the matrix $\mathfrak{B}_{A,f}$ have a nice combinatorial pattern which is explained below in cases I, II and III. 
\end{itemize}
\noindent \underline{Case I:} For the rows corresponding to elements of the type $A\frac{[-3]^{m+x}}{[1]^{m}[2]^{x}}$ of $\mathcal{R}(A,f)$ where $m$ and $x$ are integers $\geq 1$. In such a row,\\
the entry in the column corresponding to the element $A\frac{[-3]^{a+b}}{[1]^{a}[2]^{b}}$ of $P_{A,f}$ is 
\begin{center}
$\ ^{b-x}C_{m-a}$ if $0\leq a\leq m$ and $b\geq x$\\
and $0$ otherwise.
\end{center}
In the same row, the entry in the column corresponding to the element $A\frac{[-2]^{a}[-3]^{b}}{[1]^{a+b}}$ of $P_{A,f}$ (where $a$ and $b$ are integers $\geq 0$) is 
\begin{center}
$0$ if $m\geq M_{A}(\frac{-3}{1})$,\\
$(-1)^{(a+b-m)}\ ^{(a+b-m-1)}C_{a}$ if $b\neq 0$ and $m<M_{A}(\frac{-3}{1})$\\
and $0$ if $b=0$ and $m<M_{A}(\frac{-3}{1})$.
\end{center}
where $M_{A}(\frac{-3}{1})$ is the least positive integer $M$ for which $A[\frac{-3}{1}]^M$ is convergent and $A[\frac{-3}{1}]^{M-1}(-3)$ has no negative powers. When $b\neq 0$ and $m<M_{A}(\frac{-3}{1})$, $x$ has to be equal to $1$.\\
In the same row, the entry in the column corresponding to the element $A\frac{[-2]^{a+b}}{[1]^{a}[3]^{b}}$ (where $a\geq 1$ and $b\geq 0$ are integers) of $P_{A,f}$ is $0$.

\noindent \underline{Case II:} For the rows corresponding to elements of the type $A\frac{[-2]^{m+x}}{[1]^{m}[3]^{x}}$ of $\mathcal{R}(A,f)$ where $m$ and $x$ are integers $\geq 1$. In such a row,\\
the entry in the column corresponding to the element $A\frac{[-2]^{a+b}}{[1]^{a}[3]^{b}}$ (where $a\geq 1$ and $b\geq 0$ are integers) of $P_{A,f}$ is 
\begin{center}
$\ ^{b-x}C_{m-a}$ if $0\leq a\leq m$ and $b\geq x$\\
and $0$ otherwise.
\end{center}
And the entry in the column corresponding to the element $A\frac{[-3]^{a}[-2]^{b}}{[1]^{a+b}}$ of $P_{A,f}$ (where $a$ and $b$ are integers $\geq 0$) is 
\begin{center}
$0$ if $m\geq M_{A}(\frac{-2}{1})$,\\
$(-1)^{(a+b-m)}\ ^{(a+b-m-1)}C_{a}$ if $b\neq 0$ and $m<M_{A}(\frac{-2}{1})$\\
and $0$ if $b=0$ and $m<M_{A}(\frac{-2}{1})$.
\end{center}
where $M_{A}(\frac{-2}{1})$ is the least positive integer $M$ for which $A[\frac{-2}{1}]^M$ is convergent and $A[\frac{-2}{1}]^{M-1}(-2)$ has no negative powers. When $b\neq 0$ and $m<M_{A}(\frac{-2}{1})$, $x$ has to be equal to $1$.\\
In the same row, the entry in the column corresponding to the element $A\frac{[-3]^{a+b}}{[1]^{a}[2]^{b}}$ (where $a\geq 1$ and $b\geq 0$ are integers) of $P_{A,f}$ is $0$.

\noindent\underline{Case III:} For the row corresponding to the element $A\frac{[-2]}{[1]}$.\\
The entry in the column corresponding to the element $A\frac{[-2]^{a+b}}{[1]^{a}[3]^{b}}$ (where $a$ and $b$ are integers $\geq 1$) of $P_{A,f}$ is $0$. The entry in the column corresponding to the element $A\frac{[-3]^{a+b}}{[1]^{a}[2]^{b}}$ (where $a\geq 1$ and $b\geq 1$ are integers) of $P_{A,f}$ is also $0$.
And the entry in the column corresponding to the element $A\frac{[-2]^{a}[-3]^{b}}{[1]^{a+b}}$ of $P_{A,f}$ (where $a$ and $b$ are integers $\geq 0$) is $(-1)^{(a+b+1)}\ ^{(a+b)}C_{a}$.
\section{Irrationality is possible!}\label{s.irrationality}
$\mathfrak{B}_{A,f}\mathfrak{Y}=\mathfrak{e}$ is an infinite system of linear equations. The indexing set of the rows of the matrix $\mathfrak{B}_{A,f}$ is $\mathcal{R}(A,f)$, which consists of elements of the following types:
\begin{itemize}
 \item $A\frac{[-2]}{[1]}$.\\
 \item $A(\frac{[-3]}{[1]})^{m}\frac{[-3]}{[2]}$ where $0\leq m\leq M_{A}(\frac{-3}{1})-1$.\\
 \item For every fixed $m\geq M_{A}(\frac{-3}{1})$, elements of the type $A(\frac{[-3]}{[1]})^{m}(\frac{[-3]}{[2]})^{x}$ where $x$ is a positive integer lying between an interval $[a_m,b_m]$ where $a_m$ and $b_m$ are some positive integers related to the range in which mutators of the type $A(\frac{[-3]}{[1]})^{m}(\frac{[-3]}{[2]})^{N}[-3]$ lie. Details of the definition of the integers $a_m$ and $b_m$ are omitted. \\  
 \item $A(\frac{[-2]}{[1]})^{m}\frac{[-2]}{[3]}$ where $1\leq m\leq M_{A}(\frac{-2}{1})-1$.\\
 \item For every fixed $m\geq M_{A}(\frac{-2}{1})$, elements of the type $A(\frac{[-2]}{[1]})^{m}(\frac{[-2]}{[3]})^{y}$ where $y$ is a positive integer lying between an interval $[p_m,q_m]$ where $p_m$ and $q_m$ are some positive integers related to the range in which mutators of the type $A(\frac{[-2]}{[1]})^{m}(\frac{[-2]}{[3]})^{N}[-2]$ lie. Details of the definition of the integers $p_m$ and $q_m$ are omitted. 
\end{itemize}
Given any monomial $A\in\mathfrak{M}$, we have this infinite set 
\begin{center}
$$\{M_{A}(\frac{-3}{1}), M_{A}(\frac{-2}{1})\}\cup\{a_m|m\geq M_{A}(\frac{-3}{1})\}\cup\{b_m|m\geq M_{A}(\frac{-3}{1})\}$$\\
$$\cup\{p_m|m\geq M_{A}(\frac{-2}{1})\}\cup\{q_m|m\geq M_{A}(\frac{-2}{1})\}$$
\end{center}
of integers corresponding to the system $\mathfrak{B}_{A,f}\mathfrak{Y}=\mathfrak{e}$ of linear equations. Let us denote this set by $\mathbb{Z}_{A}$. But there may exist many monomials $A$ in the set $\mathfrak{M}$ for which set $\mathbb{Z}_{A}$ is the same. 

Let $\mathbb{Z}_{\mathfrak{M}}$ denote the set $\{\mathbb{Z}_{A}|A\in\mathfrak{M}\}$. Fix an arbitrary element $D$ of $\mathbb{Z}_{\mathfrak{M}}$. There will be many systems $\mathfrak{B}_{A,f}\mathfrak{Y}=\mathfrak{e}$ of linear equations corresponding to every monomial $A$ associated to $D$. For computing the Hilbert-Kunz function corresponding to a given trinomial hyper surface, we need to determine that given any such set $D$, how many systems corresponding to it are solvable and how many are not. The number of systems which are not solvable will contribute to a positive count ($+1$) in the formula for the Hilbert-Kunz function. 

For determining solvability of any such system, we can break it into infinitely many sub-systems corresponding to the following group of row indices:\\
(i) The group formed by the row indices $A\frac{[-2]}{[1]}$ and $A(\frac{[-3]}{[1]})^{m}\frac{[-3]}{[2]}$ where $0\leq m\leq M_{A}(\frac{-3}{1})-1$.\\
(ii) For each fixed positive integer $m\geq M_{A}(\frac{-3}{1})$, the group formed by the row indices $A(\frac{[-3]}{[1]})^{m}(\frac{[-3]}{[2]})^{x}$ where $x$ is a positive integer lying between the interval $[a_m,b_m]$.\\
(iii) The group formed by the row indices $A(\frac{[-2]}{[1]})^{m}\frac{[-2]}{[3]}$ where $1\leq m\leq M_{A}(\frac{-2}{1})-1$.\\
(iv) For each fixed positive integer $m\geq M_{A}(\frac{-2}{1})$, the group formed by the row indices $A(\frac{[-2]}{[1]})^{m}(\frac{[-2]}{[3]})^{y}$ where $y$ is a positive integer lying between the interval $[p_m,q_m]$.\\
The solvability of each of these infinitely many sub-systems can be determined using techniques similar to the case of disjoint-term trinomial hyper surfaces (as mentioned in \cite{Su1} and \cite{Su2}). So for the fixed element $D$ of $\mathbb{Z}_{\mathfrak{M}}$, the set of all monomials $A$ associated to it for which the corresponding systems $\mathfrak{B}_{A,f}\mathfrak{Y}=\mathfrak{e}$ are solvable will be equal to the intersection of the sets of all monomials $A$ (corresponding to $D$) for which each of the infinitely many sub-systems (of the $4$ types as mentioned above) is solvable. 

Hence for the given fixed element $D$ of $\mathbb{Z}_{\mathfrak{M}}$, the number of monomials $A$ associated to it for which the corresponding systems are not solvable will equal a limit (as $m$ tends to infinity) of a sequence of rational numbers, where every term in this rational number sequence corresponds to the number of monomials $A$ (associated to $D$) for which the subsystem associated to some natural number $m$ and some group of row indices (of type (i), (ii), (iii) or (iv)) is not solvable. This limit can be an irrational number in most of the situations.

We know that if $d$ is the dimension of the trinomial hyper surface under consideration (over any field of prime characteristic $p>0$), then the corresponding Hilbert-Kunz multiplicity is the coefficient of $p^{nd}$ in the formula of the Hilbert-Kunz function. This coefficient can become irrational due to the reason(s) mentioned in the previous paragraph. This leads to the possibility that there can be examples of trinomial hyper surfaces for which the corresponding HK multiplicity is irrational.
\providecommand{\bysame}{\leavevmode\hbox
to3em{\hrulefill}\thinspace}
\providecommand{\MR}{\relax\ifhmode\unskip\space\fi MR }
\providecommand{\MRhref}[2]{%
  \href{http://www.ams.org/mathscinet-getitem?mr=#1}{#2}
} \providecommand{\href}[2]{#2}

\end{document}